\documentclass[12pt]{amsart}

\usepackage{color}
\usepackage{amstext}
\usepackage{amsthm}
\usepackage{amsrefs}
\usepackage{amsmath}
\usepackage{amssymb}
\usepackage{latexsym}
\usepackage{amsfonts}

\begin{document}
%Commands Used%

\newcommand{\ci}[1]{_{ {}_{\scriptstyle #1}}}

\newcommand{\norm}[1]{\ensuremath{\|#1\|}}
\newcommand{\abs}[1]{\ensuremath{\vert#1\vert}}
\newcommand{\p}{\ensuremath{\partial}}
\newcommand{\pr}{\mathcal{P}}

\newcommand{\pbar}{\ensuremath{\bar{\partial}}}
\newcommand{\db}{\overline\partial}
\newcommand{\D}{\mathbb{D}}
\newcommand{\B}{\mathbb{B}}
\newcommand{\Sp}{\mathbb{S}}
\newcommand{\T}{\mathbb{T}}
\newcommand{\R}{\mathbb{R}}
\newcommand{\Z}{\mathbb{Z}}
\newcommand{\C}{\mathbb{C}}
\newcommand{\N}{\mathbb{N}}
\newcommand{\F}{\mathcal{F}}
\newcommand{\scrH}{\mathcal{H}}
\newcommand{\scrL}{\mathcal{L}}
\newcommand{\td}{\widetilde\Delta}

\newcommand{\La}{\langle }
\newcommand{\Ra}{\rangle }
\newcommand{\rk}{\operatorname{rk}}
\newcommand{\card}{\operatorname{card}}
\newcommand{\ran}{\operatorname{Ran}}
\newcommand{\im}{\operatorname{Im}}
\newcommand{\re}{\operatorname{Re}}
\newcommand{\tr}{\operatorname{tr}}
\newcommand{\vf}{\varphi}
\newcommand{\f}[2]{\ensuremath{\frac{#1}{#2}}}

\def\co#1{\textcolor{red}{#1}}
\def\gr#1{\textcolor{green}{#1}}
\def\ov{\overline}
\def\BP{Blaschke product}
\def\CNBP{Carleson-Newman Blaschke product}
\def\IBP{interpolating Blaschke product}
\def\e{\varepsilon}

%%%%%%%%%%%%%%%%%%%%%%%%%%%%

\numberwithin{equation}{section}

\newtheorem{thm}{Theorem}[section]
\newtheorem{lm}[thm]{Lemma}
\newtheorem{cor}[thm]{Corollary}
\newtheorem{conj}[thm]{Conjecture}
\newtheorem{prob}[thm]{Problem}
\newtheorem{prop}[thm]{Proposition}
\newtheorem{obs}[thm]{Observation}
\newtheorem{quest}[thm]{Question}
\newtheorem*{prop*}{Proposition}

\theoremstyle{remark}
\newtheorem{rem}[thm]{Remark}
\newtheorem*{rem*}{Remark}

\title[Stable Ranks of $H^\infty_\R$ and $A_\R$]{The Bass and Topological Stable Ranks of $H^\infty_\R(\D)$ and $A_\R(\D)$}

\author[R. Mortini]{Raymond Mortini}
\address{Raymond Mortini, D\'{e}partement de Math\'{e}matiques\\ LMAM, UMR 7122, Universit\'{e} Paul Verlaine\\ Ile du Saulcy\\  F-57045 Metz, France}
\email{mortini@univ-metz.fr}

\author[B. D. Wick]{Brett D. Wick$^*$}
\address{Brett D. Wick, Department of Mathematics\\ University of South Carolina\\ LeConte College\\ 1523 Greene Street\\ Columbia, SC USA 29208}
\address{Present:  The Fields Institute\\ 222 College Street, 2nd Floor\\ Toronto, Ontario\\ M5T 3J1 Canada}
\email{wick@math.sc.edu}
\thanks{$*$ Research supported in part by National Science Foundation Grant DMS -- 0752703.}


\begin{abstract}
In this note we prove that the Bass stable rank of $H^\infty_\R(\D)$ is two.  This establishes the validity of a conjecture by S. Treil.  We accomplish this in two different ways, one by giving a direct proof, and the other, by first showing that the topological stable rank of $H^\infty_\R(\D)$ is two.  
 We apply these results to give new proofs of results by R. Rupp and A. Sasane stating that the Bass stable rank of $A_\R(\D)$ is two and the topological stable rank of $A_\R(\D)$ is two, settling a conjecture by the second author.  We also present a $\bar\partial$-free proof of the second author's characterization of the reducible pairs in $A_\R(\D)$.
\end{abstract}

\maketitle

\section{Introduction and Main Results}
\label{Intro.}
The stable rank of an algebra (also called the Bass stable rank) was introduced by H. Bass in \cite{Bass} to assist in computations of algebraic K-Theory.  We recall the definition of the stable rank.

Let $A$ be a commutative topological algebra with  unit $e$. 
An $n$-tuple $a=(a_1,\cdots,a_n)\in A^n$ is called \textit{unimodular} or \textit{invertible} if there exists an $n$-tuple $b\in A^n$ such that $\sum_{j=1}^{n}a_jb_j=e$.  The set of all invertible $n$-tuples in $A^n$ is denoted by $U_n(A)$.
An $(n+1)$-tuple $a$ is called \textit{stable} or \textit{reducible} if there exists an $n$-tuple $x$ such that the $n$-tuple $(a_1+x_1a_{n+1},\ldots, a_{n}+x_{n}a_{n+1})$ is invertible.  The \textit{stable rank} (also called \textit{bsr}(A) in the literature) of the algebra $A$ is the least integer $n$ such that every invertible $(n+1)$-tuple is reducible.  

Another invariant of the algebra $A$ that we will be interested in is the \textit{topological stable rank} (frequently called \textit{tsr}(A) in the literature). Recall that the topological stable rank of  $A$  is the least integer $n$ such that $U_n(A)$ is dense in $A^n$. We always have that \textit{bsr}(A) $\leq$ \textit{tsr}(A), see \cite{CL}.

These stable ranks are purely algebraic concepts but can be combined with analysis when studying commutative Banach algebras of functions.  In this context, the stable rank is related to the zero sets of ideals, and the spectrum of the Banach algebra.  The stable ranks for different algebras of analytic functions have been considered by many authors.  Jones, Marshall and Wolff, \cite{JMW}, showed that $\text{bsr}(A(\D))=1$.  The computation of stable ranks for sub-algebras of the disc algebra $A(\D)$  was done by Corach--Su\'arez, \cite{CS}, and Rupp \cite{R}.  The determination of the topological stable rank of $A(\D)$ can be found in the paper \cite{JMW} and can easily be seen from elementary results in complex function theory.

For the Banach algebra of all bounded analytic functions in the open unit 
disc $\D$, $H^\infty(\D)$, the classification of its invertible $n$-tuples and its stable rank are well understood.  Carleson's Corona Theorem, see \cite{Carleson}, can be phrased as an $n$-tuple $(f_1,\ldots,f_n)\in H^\infty(\D)^n$ is invertible if and only if it satisfies the Corona condition,
$$
\inf_{z\in\D}\left(\abs{f_1(z)}+\cdots+\abs{f_n(z)}\right)=\delta>0.
$$

The stable rank of $H^\infty(\D)$ was computed by S. Treil and is one of the motivations for this paper.  Treil's result is the following theorem:

\begin{thm}[S. Treil \cite{TreilStable}]
\label{stableHinfty}
Let $f_1,f_2\in H^\infty(\D)$ be such that $||f_j||_\infty\leq 1$ and $\inf_{z\in\D}(\abs{f_1(z)}+\abs{f_2(z)})=\delta>0$.  Then there exist $g_1, g_2, g_1^{-1}\in H^\infty(\D)$ with $\norm{g_1}_\infty,\norm{g_2}_{\infty}$ and $\norm{g_1^{-1}}_\infty$ controlled by $C(\delta)$, a constant only depending on $\delta$, and  
$$
1=f_1(z)g_1(z)+f_2(z)g_2(z),\quad\forall z\in\D.
$$
\end{thm}

The topological stable rank of $H^\infty(\D)$ was determined by D. Su\'arez, and serves as another motivation for this paper.  The result of interest is the following,

\begin{thm}[D. Su\'arez \cite{su}]
The topological stable rank of $H^\infty(\D)$ is 2.
\end{thm}

It is apparent that Theorem \ref{stableHinfty} implies the stable rank of $H^\infty(\D)$ is one.  Questions about the stable rank of some sub-algebras of $H^\infty(\D)$ have been studied by the first author, \cite{Mortini}.  

It is possible to phrase Treil's result \cite{TreilStable} in the language of Control Theory.  In this language, the result can be viewed as saying that it is possible to stabilize (\textit{in the sense given above}) a linear system (\textit{the Corona data, viewed as a rational function}) via a stable (\textit{analytic}) controller.  But, in applications of Control Theory, the linear systems and transfer functions  have real coefficients, so in this context Treil's result is physically meaningless.  From the point of view of Control Theory, it is important to know whether results like Theorem \ref{stableHinfty} hold, but for a more physically meaningful algebra. This serves as further motivation for our paper.
We are interested in questions related to the stable rank of a natural sub-ring of $H^\infty(\D)$, the real Banach algebra $H^\infty_\R(\D)$.  In particular, we want to know whether  some variant of Theorem \ref{stableHinfty} holds for this algebra. 

First, recall that $H^\infty_\R(\D)$ is the subset of $H^\infty(\D)$ with the additional property that the Fourier coefficients of an element of $H^\infty_\R(\D)$ must be real.  This property can be captured by the following symmetry condition
$$
f(z)=\overline{f(\overline{z})}\quad\forall z\in\D.
$$
This condition implies that the functions in $H^\infty_\R(\D)$ possess a symmetry that is not present  for general $H^\infty(\D)$ functions.

Carleson's Corona result is inherited by the algebra $H^\infty_\R(\D)$.  More precisely, it is an immediate application of the usual Corona Theorem and the symmetry properties of $H^\infty_\R(\D)$ to show that an 
$n$-tuple $(f_1,\ldots,f_n)\in H^\infty_\R(\D)^n$ is invertible if and only if it satisfies the Corona condition,
$$
\inf_{z\in\D}\left(\abs{f_1(z)}+\cdots+\abs{f_n(z)}\right)=\delta>0.
$$
Indeed, one direction is immediate, and in the other direction, if we know that
$$
\inf_{z\in\D}\left(\abs{f_1(z)}+\cdots+\abs{f_n(z)}\right)=\delta>0,
$$
then we can find a solution $(g_1,\ldots,g_n)\in H^\infty(\D)^n$.  We then symmetrize the $g_j$ via the operation
$$
g_j^{\#}(z):=\f{g_j(z)+\overline{g_j(\overline{z})}}{2}.
$$
Then $g_j^{\#}\in H^\infty_\R(\D)$ and $(g^{\#}_1,\dots,  g^{\#}_n)\in H^\infty_\R(\D)^n$ is a Corona solution of the Bezout equation $$\sum_{j=1}^n g_j^{\#}f_j=1.$$

The question of whether the algebra $H^\infty_\R(\D)$ has Bass stable rank one was addressed by the second author in \cite{Wick1}.  

It is very easy to see that there is an additional necessary condition that must be satisfied by the invertible pairs in $H^\infty_\R(\D)$.  
Suppose that for $f_1,f_2\in H^\infty_\R(\D)$ there exist functions 
$g_1,g_2\in H^\infty_\R(\D)$ satisfying the Bezout equation $f_1g_1+f_2g_2=1$.
 Then we shall see that the real zeros of $f_1$ and $f_2$ must intertwine correctly.  Indeed, let $\lambda_1$ and $\lambda_2$ be real zeros of $f_2$.  Then we have
\begin{eqnarray*}
f_1(\lambda_1)g_1(\lambda_1) & = & 1\\
f_1(\lambda_2)g_1(\lambda_2) & = & 1.
\end{eqnarray*}
Now $f_1(\lambda_1)$ and $f_1(\lambda_2)$ must have the same sign at these zeros.  If this were not true, then without loss of generality suppose that $f_1(\lambda_1)>0>f_1(\lambda_2)$.  Then $g_1(\lambda_1)>0>g_1(\lambda_2)$.  By continuity there will exist a point $\lambda_{12}$, between $\lambda_1$ and $\lambda_2$, with $g_1(\lambda_{12})=0$.  But this contradicts the fact that $g_1^{-1}\in H^\infty_\R(\D)$.  So $f_1$ must have the same sign at real zeros of $f_2$.  We will say that $f_1$ is \textit {of constant sign} on the real zeros of $f_2$ if $f_1$ has the same sign at all real zeros of $f_2$.

This is also an intertwining condition of the zeros of $f_1$ and $f_2$.  More precisely, the function $f_1$ is of constant sign on the real zeros of $f_2$ if and only if, between two real zeros of $f_2$, there must be an even number of real zeros of $f_1$, or none.  This is called the \textit{parity interlacing property} found in Control Theory for the stabilization of a linear system.

The final motivation for this paper is the following theorem of the second author.

\begin{thm}[Wick, \cite{Wick1}]
\label{Hinfty}
Suppose $f_1,f_2\in H^\infty_\R(\D)$, $\norm{f_1}_\infty,\norm{f_2}_{\infty}\leq 1$, $f_1$ is of constant sign on the real zeros of $f_2$ and
$$
\inf_{z\in\D}\left(\abs{f_1(z)}+\abs{f_2(z)}\right)=\delta>0.
$$
Then there exists $g_1, g_1^{-1}, g_2\in H^\infty_\R(\D)$ with $\norm{g_1}_\infty,\norm{g_2}_\infty,\norm{g_1^{-1}}_\infty\leq C(\delta)$ and
$$
f_1(z)g_1(z)+f_2(z)g_2(z)=1\quad\forall z\in\D.
$$
\end{thm}

The reasoning leading to this theorem in turn implies that the stable rank of $H^\infty_\R(\D)$ is at least two because the additional condition, $f_1$ is of constant sign on the real zeros of $f_2$, is necessary to find Corona solutions with one of them invertible.  

Thus, a natural question, first conjectured by S. Treil and further popularized by the second author, \cite{Wick1}, is that the Bass stable rank of $H^\infty_\R(\D)$ is two.  

\subsection*{Main Results}
This leads to the main results of this paper which proves that the conjecture of Treil
is true.  We will establish this conjecture in two different ways.  One will be to demonstrate that the Bass stable rank of $H^\infty_\R(\D)$ is two, and the other will be to demonstrate that the topological stable rank of $H^\infty_\R(\D)$ is two.

\begin{thm}[Main Result]
\label{str}
The stable rank of $H^\infty_\R(\D)$ is 2.
\end{thm}

The other main result of this paper is the computation of the topological stable rank of $H^\infty_\R(\D)$.

 \begin{thm}[Main Result]
 \label{tsr}
 The topological stable rank $H^\infty_\R(\D)$ is 2.
 \end{thm}

We remark at this point that Theorem \ref{tsr} implies Theorem \ref{str} by well known facts.  However, we will provide a direct proof of Theorem \ref{str}.

The other results of this paper address questions related to the algebra $A_\R(\D)$, the real disc algebra.  Recall that the disc algebra is the set of bounded analytic functions that admit a continuous extension to the boundary of the disc.  The real disc algebra is the collection of functions in $A(\D)$ that satisfy the additional symmetry property of being real on the real axis.

For this algebra, we give several results.  Using Theorem \ref{str} and Theorem \ref{tsr} we compute the  Bass and topological stable ranks of $A_\R(\D)$.  These results were conjectured by the second author in \cite{Wick2} and then proven by Rupp and Sasane in \cite{RS1}.  See also \cite{RS2}.
  
\begin{thm}[Rupp and Sasane, \cite{RS1}]
\label{rsrank}
The Bass stable rank of $A_\R(\D)$ is 2.
\end{thm}

\begin{thm}[Rupp and Sasane, \cite{RS1}]
\label{rstrank}
The topological stable rank of $A_\R(\D)$ is 2.
\end{thm}

Again, we remark that Theorem \ref{rstrank} implies Theorem \ref{rsrank} by well known results.  However, we will deduce both theorems directly.
  
Finally, we give a $\db$-free proof of the main result in \cite{Wick2} by the second author.  Namely, we prove that

\begin{thm} [Wick, \cite{Wick2}] 
\label{wick}
An invertible pair $(f,g)$ in  $A_\R(\D)^2$ is reducible if and only if  $f$ has constant sign on the real zeros of $g$, or, which is equivalent, if and only if between two real zeros of $g$ in $\overline\D$ there is
an even number of zeros of $f$ (or none).
\end{thm}

\section{The Bass and Topological Stable Ranks of $H^\infty_\R(\D)$}

We first collect some key propositions that will be used in the proof of the main result.  The key idea is that the functions in $H^\infty_\R(\D)$ possess a certain symmetry that can be exploited.  Using this symmetry, it is possible to take arbitrary functions from $H^\infty(\D)$ and construct functions that are in $H^\infty_\R(\D)$.  The main point will be to do so without introducing ``extra zeros'' into the function.

We first state a lemma that applies, in general, to any commutative Banach algebra $A$.  See for example, \cite[p. 636]{CS} and \cite[p. 608]{cs2}.

\begin{lm}
\label{clopen}
 Let $A$ be  a commutative unital Banach algebra. Then, for $g\in A$, the set 
$$
\mbox{$R(g)=\{f\in A: (f,g)$ is reducible$\}$}
$$
 is open-closed inside $I(g)=\{f\in A: (f,g)$ is an invertible pair in $A^2\}$.
 In particular,
if $\phi: [0,1]\to I(g)$ is a continuous curve and $(\phi(0),g)$ is reducible, then $(\phi(1),g)$
is reducible.
\end{lm}
  
For $u\in H^\infty(\D)$ let $u^*(z):=\overline{u(\overline z)}$.  The reader should not be confused by the similarity with $u^{\#}(z):=\frac{u(z)+u^*(z)}{2}$.
    
 \begin{prop}\label{square}
 Let $(u,g)$ be an invertible pair in $H^\infty_\R(\D)^2$. Then $(u^2,g)$ is reducible in $H^\infty_\R(\D)$.
 \end{prop}
 
\begin{proof}
First, we note that the necessary condition for being reducible in  $H^\infty_\R(\D)$, namely
that $u^2$ is of constant sign on the real zeros of $g$, is trivially satisfied.

Since $H^\infty(\D)$ has Bass stable rank one, there exists $h\in H^\infty(\D)$ so that $v:=u+hg$ is invertible 
in $H^\infty(\D)$. Hence $v^*=u+h^*g$ is invertible in $H^\infty(\D)$. Multiplying both equations, yields
that 
$$
u^2+ g\bigl( u(h+h^*)+hh^*g\bigr) = vv^*
$$
is invertible in $H^\infty(\D)$. But, the factors $h+h^*$ and $hh^*$ belong to  $H^\infty_\R(\D)$.
Hence $(u^2,g)$ is reducible in $H^\infty_\R(\D)$.
\end{proof}

\begin{cor}
Let $(u,g)$ be an invertible pair in $H^\infty_\R(\D)^2$. Suppose that $u$ belongs to the closure
of the invertibles in $H^\infty_\R(\D)$. Then $(u,g)$ is reducible in $H^\infty_\R(\D)$.
\end{cor}

\begin{proof}
This immediately follows from the representation above:
$$
u\cdot u +g F\in(H^\infty_\R(\D))^{-1},
$$
and the fact that $u$ is uniformly close to some $v\in (H^\infty_\R(\D))^{-1}$.
\end{proof}

\begin{cor}\label{zerofree}

Let $(f,g)$ be an invertible pair in $H^\infty_\R(\D)^2$ and suppose that $f$ has no zeros in $\D$.
Then  $(f,g)$ is reducible in $H^\infty_\R(\D)$.
\end{cor}

\begin{proof}
Without loss of generality, we may suppose that  $f$ is strictly positive on $]-1,1[$. Then the assertion follows from Proposition \ref{square} above, by noticing that $f$ has a square root
in  $H^\infty_\R(\D)$:  
$$
\sqrt{f(z)}=\exp\frac{1}{2}\left[  \int_0^{2\pi} \frac{e^{i t}+z}{e^{it}-z}\biggl(\log|f(e^{it})| \frac{dt}{2\pi}- d\mu_s(t)\biggr)\right]
$$
where $\mu_s$ is a positive singular measure with $\mu_s(E)=\mu_s(\overline E)$ for any Borel set $E$
on $\T$.
\end{proof}

Before the next proposition, we need the definition of the zero set of a function.  Given a holomorphic function $B$, we define
$$
Z_\D(B):=\{z\in\D: B(z)=0\}.
$$
Also, if $A=H^\infty(\D)$ or $A=A(\D)$, then the zero set of the function $f\in A$ 
on the spectrum $M(A)$ of $A$ is defined as

$$
Z(f):=\{m\in M(A): \hat f(m)=0\},
$$
where $\hat f: m\mapsto m(f)$, $m\in M(A)$,  denotes the Gelfand transform of $f$.
In the sequel we shall always identify $f$ with $\hat f$.

\begin{prop}

Let $(f,g)$ be an invertible pair in $H^\infty_\R(\D)^2$. Suppose that  $f>0$ on $Z(g)$ in the spectrum of 
$H^\infty(\D)$. Then $(f,g)$ is reducible in $H^\infty_\R(\D)$.
\end{prop}

\begin{proof}
We may assume that $||f||_\infty<\frac{1}{2}$. Consider the continuous curve $[0,1]\to  H^\infty_\R(\D), 
\epsilon\mapsto f+\epsilon$. Since $f>0$ on $Z(g)$ we also have that $f+\epsilon>0$ on $Z(g)$. So
$(f+\epsilon,g)$ is an invertible pair for all $\epsilon\in [0,1]$. Now $f+1$ itself is invertible in
$ H^\infty_\R(\D)$; hence $(f+1,g)$ is reducible. By Lemma \ref{clopen},  $(f+\epsilon,g)$ is reducible for every $\epsilon$;
in particular $(f,g)$. 
\end{proof}

\begin{prop}\label{real}
Let $(f,g)$ be an invertible pair in $H^\infty_\R(\D)^2$. Suppose that $f$ has no real zeros. Then
$(f,g)$ is reducible in $H^\infty_\R(\D)$.
\end{prop}

\begin{proof}
Write $f=BF$ where $F$ is zero-free and $B$ is the \BP\ associated with the zeros of $f$.
  Since the zeros of $B$ are symmetric with respect to the real axis, $B$ and $F$ are in $H^\infty_\R(\D)$. Thus $(B,g)$ and $(F,g)$ are invertible pairs
in $H^\infty_\R(\D)^2$.  Without loss of generality, we may suppose that  
$F>0$ on $]-1,1[$.
By Corollary \ref{zerofree}, $(F,g)$ is reducible.  It remains to show
that $(B,g)$ is reducible. Let $C$ be the \BP\ formed with the zeros of $B$
on the upper-half plane and let $C^*$ be defined as $C^*(z):=\overline{C(\overline z)}$. Then
$B=CC^*$.  Since $(C,g)$ is an invertible pair in $H^\infty(\D)^2$, there exists by Treil's Theorem, see \cite{TreilStable}, an $h\in H^\infty(\D)$ and $v\in (H^\infty(\D))^{-1}$ such that $C+hg =v$. But $v^*= C^*+h^*g$.
Multiplying both equations yields
$$
CC^*+\bigl((Ch^*+C^*h)+hh^*g\bigr)g  =vv^*
$$
Thus $(B,g)$ is reducible in $H^\infty_\R(\D)$.
\end{proof}

\begin{prop}
\label{realzeros} 
Let $B$ be a \BP\ having only  real zeros, all of them being simple, and let $f\in H^\infty_\R(\D)$ 
  be such that  $(f,B)$ is  an invertible pair in  $H^\infty_\R(\D)^2$. Then there exist
  two factors $B_1$ and $B_2$ of $B $ such that $B=B_1B_2$ and such that
  $(f,B_j)$ are reducible in $H^\infty_\R(\D)$ for $j=1,2$. Moreover,  $B_1$ can be chosen to be 
  the sub-product of $B$ for which $f>0$ on $Z_\D(B_1)$ and $B_2$ the factor
   for which $f<0$ on $Z_\D(B_2)$.
\end{prop}

\begin{proof}
Let $B_1$ and $B_2$ be as above. By \cite{TreilStable}, there exists $F\in(H^\infty(\D))^{-1}$ and $h\in H^\infty(\D)$ such that
$Ff+hB=1$. Hence $F^*f+h^*B=1$. Multiplying both equalities, shows that $FF^*f^2\equiv 1$
on $Z_\D(B)$. Now $FF^*$ is an invertible function in $H^\infty_\R(\D)$, hence, as above, has
 a square root in $H^\infty_\R(\D)$, say $v^2=FF^*$. 
 We may assume that $v(0)>0$. Then $(vf-1)(vf+1)=FF^*f^2-1\equiv 0$
 on $Z_\D(B)$.  We deduce that $vf-1 =b_1q_1$ and that $vf+1=b_2q_2$ for some
 factors $b_j$ of $B$ with $b_1b_2=B$  and some $q_j\in H^\infty(\D)$. Since $vf\pm 1\in H^\infty_\R(\D)$, 
  $q_j\in H^\infty_\R(\D)$, too. Now $vf =(-1)^{j-1} $ on $Z_\D(b_j)$. Since $v$ is invertible in $H^\infty_\R(\D)$,
for each $j$, $f$ has constant sign on $Z_\D(b_j)$.  Thus $b_j=B_j$. 
\end{proof}

An immediate, and important, Corollary of the above result is the following statement.

\begin{prop}
\label{realzeros2}
Let $B$ be a \BP\ having only simple, real zeros and let $f\in H^\infty_\R(\D)$ 
  be such that  $(f,B)$ is  an invertible pair in  $H^\infty_\R(\D)^2$.  Suppose that $f$ has constant sign on $Z_\D(B)$, then $(f,B)$ is reducible.
\end{prop}

Let $\F$ be the set of finite products of \IBP s.  These functions are also often called \CNBP s.
The next Lemma is well known for people familiar with Hoffman's theory
on the maximal ideal space. 
For the reader's convenience, we present a proof.

\begin{lm}\label{cnbp}
Let $B\in\F\cap H^\infty_\R(\D)$ and let $f\in H^\infty_\R(\D)$. Then for every $\e>0$ small enough,
$B+\e f$ writes as $C v$, where $C\in\F\cap H^\infty_\R(\D)$, and where $v$ is an invertible
outer function in  $H^\infty_\R(\D)$.
\end{lm}

\begin{proof}
We shall use a result of Guillory, Izuchi and Sarason, see \cite{gis}, which tells us that an inner function $B$ is in $\F$ if and only if $B$ does not vanish on the compact set $X$ of trivial points
in the maximal ideal space $M(H^\infty(\D))$ of $H^\infty(\D)$.  So, $B\in \F$ implies that 
$$|B+\e f|\geq \eta-\e||f||_\infty=:\delta>0$$
 for small $\e>0$ on $X$.  Let $Cv$ be the inner-outer factorization of $B+\e f$. Since the inner function $C$ does not vanish on $X$, it belongs to 
 $\F$ by \cite{gis} (In particular $C$ has no singular inner factors).
Since $Cv\in H^\infty_\R(\D)$, the  zeros of $C$ are symmetric; hence  $C$
 and $v$ are in $H^\infty_\R(\D)$. Moreover, since $|C|=1$ on the Shilov boundary
 $\partial H^\infty(\D)$
 of $H^\infty(\D)$, a set that can be identified with the spectrum of $L^\infty(\T)$ and that is properly contained in $X$, we see that $|v|\geq \delta$ on 
 $\partial H^\infty(\D)$. Thus $v$ is invertible in $L^\infty(\T)$ and hence (being outer)
 $v$ is invertible in $H^\infty(\D)$. 
\end{proof}

The following theorem, which  follows from \cite[Corollary 4, p.~290]{mn} and \cite[Theorem 1, p.~288]{mn} will be used rather frequently in the course of the proof of our main results.

\begin{thm}[Mortini, Nicolau, \cite{mn}]\label{frost}
Any Blaschke product whose zeros belong to a finite union of   Stolz angles belongs to the set of inner functions all of whose non-zero Frostman shifts are \CNBP s.
\end{thm}

The following is an analogue of Su\'arez's  Lemma \cite[Lemma 4.3]{su} for
$H^\infty_\R(\D)$. There are functions
$B\in \F\cap H^\infty_\R(\D)$ that cannot be decomposed into a finite
product
of \IBP s themselves being symmetric.  A simple example is obtained by taking an interpolating
sequence $\{a_n\}_{n\in\N}$
in the upper half disc converging tangentially to 1 so that
$\rho(a_n,\overline{a_n})\to 0$, with $\rho(z,w):=\left|\frac{z-w}{1-\ov z w}\right|$ the pseudohyperbolic metric, and by considering the \BP\ $B$ with zeros $\{a_n,\overline{a_n}: n\in \N\}$.  Thus the next lemma does not follow directly from \cite[Lemma 4.3]{su} and needs a somewhat tricky
argument to prove.

\begin{lm}\label{sua}
Let $B\in\F$ and suppose that  $B\in H^\infty_\R(\D)$. Moreover,
let $\varphi\in H^\infty_\R(\D)$. Then, for every $\epsilon>0$ there exists
$\varphi_1\in H^\infty_\R(\D)$
and $B_1\in H^\infty_\R(\D)$  with $||\varphi-\varphi_1||_\infty<\epsilon$ and $||B-B_1||_\infty<\epsilon$
such that $(B_1,\varphi_1)$ is an invertible pair in $H^\infty_\R(\D)^2$.

\end{lm}
\begin{proof}
Let $\varphi=(cCC^*)G$ be the Riesz factorization of $\varphi$, where $c$ is the \BP\
formed with the real zeros of $\varphi$, $C$ the \BP\ with the zeros in the
upper half disc
and where $G\in H^\infty_\R(\D)$ is zero free in $\D$. Without loss of generality,
we may assume that $G>0$ on $]-1,1[$. Since $\sqrt
G\in H^\infty_\R(\D)$,
we may write $G$ as $G=FF^*$ (just let $F=\sqrt G$).

Let $0<r<1$. By Theorem \ref {frost}, % \cite{mn},  
 the Frostman shift $c_r:=\frac{c-r}{1-r c}$ belongs to $\F$. Note that $c_r\in
H^\infty_\R(\D)$.
Let $\varphi_0=c_rCC^*FF^*$. Then $||\varphi-\varphi_0||_\infty\leq \epsilon$ whenever $r$ is small enough.
Next we write $c_r$ as $c_r=pbb^*$,  where $p\in\F$ has only real zeros and
where $b\in\F$ is formed with the zeros of $c_r$ that lie in the upper
half disc.
Then $\varphi_0$ writes as $\varphi_0=pHH^*$ for some $H\in H^\infty(\D)$ without real zeros
 and $p\in  H^\infty_\R(\D)\cap\F$. In order to approximate the pair
$(B,\varphi_0)$ by invertible ones,
we shall deal with the pairs $(B,p)$ and
$(B,HH^*)$ separately.

\medskip

{\it First Step:}
Consider the pair $(B,p)$. We now mimic Su\'arez's inductive proof of
\cite[Lemma 4.3]{su}.
 Let $p=p_1\cdots p_N$, where the $p_j$ are \IBP s with real zeros only.

 Consider the case $N=1$.  Suppose that $p_1$ satisfies
 $$ (1-|z_n|^2) |p_1'(z_n)|\geq \delta>0,$$
 where $(z_n)$ is the zero sequence of $p_1$.
 Factorize $p_1=q_1q_2$ as follows:  The zeros of $q_1$ are those $z_n$
for which
 $|B(z_n)|\geq \epsilon^2\delta$ and the zeros of $q_2$ are those $z_n$ for which
 $|B(z_n)|<\epsilon^2\delta$.  Take $R_1=B-\epsilon q_1$. Then $R_1\in H^\infty_\R(\D)$.
 Moreover, if $\epsilon$ is small enough, by Lemma \ref{cnbp}, $R_1$ has the form $\widetilde R_1
v_1$, where
 $v_1$ is invertible in $H^\infty_\R(\D)$ and where $\widetilde R_1\in\F\cap
H^\infty_\R(\D)$.

 Then it can be seen, as in
\cite[Lemma 4.3]{su}, that $|R_1|$
 is bounded away from zero on the zeros of $p_1$.  Indeed, if $q_1(z_n)=0$ then 
 \begin{equation}
 \label{zero1}
 |R_1(z_n)|=|B(z_n)-\epsilon q_1(z_n)|=|B(z_n)|\geq\epsilon^2\delta
 \end{equation}
 with the last inequality holding by the definition of the zeros of $q_1$.  On the other hand, if $q_2(z_n)=0$ then
 \begin{equation}
 \label{zero2}
 |R_1(z_n)|\geq\epsilon|q_1(z_n)|-|B(z_n)|\geq\epsilon\delta-\epsilon^2\delta>0,
 \end{equation}
 again with the last inequalities following by the definitions of the zeros of $q_2$.  Hence, $(R_1, p_1)$ is an invertible pair in $H^\infty_\R(\D)^2$.  Moreover, $||B-R_1||_\infty<\epsilon$ by the definition of $R_1$.  

 Via induction (see \cite[Lemma 4.3]{su}), we get the existence of a function $R_N$ (that has the
 form $R_N=\widetilde R_N v_N$  for some  $\widetilde R_N\in\F\cap
H^\infty_\R(\D)$ and some $v_N$ which is invertible in $H^\infty_\R(\D)$ with
$||v_N||_\infty\leq 2$)
 such that $(R_N, p_1\cdots p_N)$ is
 an invertible pair in $H^\infty_\R(\D)$ and such that $||R_N-B||_\infty\leq
\epsilon$.  To sum up, we have approximated $(B,p)$ by the  invertible pair $(R_N,p)$ in $H^\infty_\R(\D)^2$.

\medskip

 {\it Second Step:} Consider the pair $(R_N, HH^*)$.  Or, even better, one can look at 
 $(\widetilde R_N,  v_N^{-1}HH^*)$, which we can write as $(\widetilde
R_N, KK^*)$ for some $K\in H^\infty(\D)$.

Using \cite[Lemma 4.3]{su}, and noticing that $\widetilde R_N\in\F$,
there exists $\widetilde K\in
H^\infty(\D)$ with
 $$
 ||K-\widetilde K||_\infty<\frac{\epsilon}{4||K||_\infty+2}
 $$ 
 so that $(\widetilde R_N,
\widetilde K)$ is an invertible pair in
 $H^\infty(\D)^2$. Note that the  Gelfand transform of  $\widetilde R_N$
in $M(H^\infty(\D))$
vanishes only at the closure of the zeros of $\widetilde R_N$ in $\D$. Thus,
the symmetry of $\widetilde R_N$ implies that $\widetilde K^*$ has no
zeros in common
 (on $M(H^\infty(\D))$ with
 $\widetilde R_N$, too. So $(\widetilde R_N, \widetilde K \widetilde
K^*)$ is an invertible pair in $H^\infty_\R(\D)^2$.
 Moreover, 
 $$
 ||KK^*-\widetilde K \widetilde K^*||_\infty\leq  \frac{\epsilon}{2}.
 $$ 

 {\it Third Step:} Now $||v_N \widetilde K \widetilde K^*-v_NKK^*||_\infty\leq \epsilon$.
 Note that $v_NKK^*=HH^*$. Thus
 we have shown that $(R_N, HH^*)=(v_N\widetilde R_N, HH^*)$   has been
approximated
 by the invertible pair
 $$(R_N,v_N \widetilde K \widetilde K^*)=
 (v_N \widetilde R_N, v_N \widetilde K \widetilde K^*).$$

 {\it Fourth Step:} Combining the first and third step
 we see that the invertibility of the pairs   $(R_N,v_N \widetilde K
\widetilde K^*)$ and
 $(R_N,p)$ implies the invertibility of $(R_N, pv_N \widetilde K
\widetilde K^*)$.
 Let $\varphi_1=pv_N \widetilde K \widetilde K^*$ and $B_1=R_N$.  Then  $(B_1,\varphi_1)$ is the invertible pair in $H^\infty_\R(\D)^2$ we were looking for; in fact

 $||B_1-B||_\infty<\epsilon$ and $||\varphi_1-\varphi||_\infty<\epsilon$, since
\begin{eqnarray*}
||\varphi-\varphi_1||_\infty & \leq & ||\varphi-\varphi_0||_\infty+||\varphi_0-\varphi_1||_\infty\\
 & < & \epsilon+||pHH^*-\varphi_1||_\infty\\
 & < & \epsilon+||HH^*-v_N\widetilde K\widetilde K^*||_\infty\\
 & = & \epsilon+ ||v_NKK^*-v_N\widetilde K \widetilde K^*||_\infty<2\epsilon.
\end{eqnarray*}
\end{proof}

With these preliminary results collected, we now give a proof of the main results.

\subsection*{The Bass Stable Rank of $H^\infty_\R(\D)$}

We now will show that the Bass stable rank of $H^\infty_\R(\D)$ is two.  To accomplish this, we need to demonstrate that every  invertible triple in $H^\infty_\R(\D)^3$ is reducible to an invertible pair in $H^\infty_\R(\D)^2$.  Additionally, we need to know that there exists an invertible pair in $H^\infty_\R(\D)^2$ that is not reducible in $H^\infty_\R(\D)$.  The second step will be easy, while the first will depend on the previous propositions.

\begin{thm}
The stable rank of $H^\infty_\R(\D)$ is 2.
\end{thm}

\begin{proof}

As was noticed in \cite{Wick1}, the stable rank of $H^\infty_\R(\D)$ is bigger than two.  In fact, consider
the invertible pair $(z, 1-z^2)$.  Suppose that  $v:=z+h(1-z^2)$ is invertible in $H^\infty_\R$
for some $h\in H^\infty_\R$. Then $\lim_{x\to 1} v(x)=1$ and $\lim_{x\to -1} v(x)=-1$.
Hence $v$  necessarily must have a real zero on $]-1,1[$, a contradiction. 

Now let $(f,g,h)$ be an invertible triple in $H^\infty_\R(\D)^3$; say 
$$xf+yg+th=1$$
 for some 
$(x,y,t)\in H^\infty_\R(\D)^3$. Let $R=yg+th$. Write $x=BX$, where $B$ is the Blaschke
 product formed with the real zeros of $x$ (if there are none, let $B\equiv1$).  
 Obviously $(X,R)$ is an invertible pair in $H^\infty_\R(\D)^2$. Since  $X$ has no real zeros,
 this pair is reducible by Proposition \ref{real}, say $X+aR=v$ for some $a\in H^\infty_\R(\D)$
 and $v\in (H^\infty_\R(\D))^{-1}$.  Thus
 $$
 1= B(v-aR)f +R= vBf + (1-aBf)R;   
 $$
 hence for some $p, q\in H^\infty_\R(\D)$,  the sum  $Bf+ pg+qh$ is an invertible  function in $H^\infty_\R(\D)$.
 
 By Theorem \ref{frost}, for every $r\in \;]0,1[$, the Frostman shift $B_r:=\frac{B-r}{1-r B}$ of  
 $B$ is a finite product of  \IBP s. Thus, by taking $r$ small enough, we may assume
 that $B_r f+pg+qh\in (H^\infty_\R(\D))^{-1}$. 
 
Now, we apply Lemma \ref{sua} to the function $B_r$ and $p$.  Thus, there exist functions $B_1\in H^\infty_\R(\D)$ and  $p_1\in H^\infty_\R(\D)$ with $||p - p_1||_\infty<\epsilon$ and $||B_r-B_1||_\infty<\epsilon$ so that $(B_1,p_1)$ is an invertible pair in $H^\infty_\R(\D)^2$.
 
 Thus, for $\epsilon$ small enough,  $u:=B_1 f+p_1g+qh$ is invertible in $H^\infty_\R(\D)$. Moreover,
 there exist $k_1$ and $k_2\in H^\infty_\R(\D)$ such that $q=k_1 B_1+k_2p_1$.
 Hence 
$$
u= B_1 f+p_1g +(k_1B_1+k_2p_1)h = B_1(f+k_1h) +  p_1(g +k_2h).
$$ 
 Thus $(f,g,h)$ is reducible in $H^\infty_\R(\D)$.
 \end{proof}
 
\subsection*{The Topological Stable Rank of $H^\infty_\R(\D)$}

Our next goal is to prove that the topological stable rank of
$H^\infty_\R(\D)$ is two.
 Recall that the topological stable rank of a unital commutative Banach
algebra $A$ is the smallest integer $n$ for which the set of invertible $n$-tuples
is dense in $A^n$.

We will let $U_2(H^\infty_\R(\D))$ denote the set of all invertible pairs in $H^\infty_\R(\D)^2$.

 \begin{thm}
 The topological stable rank $H^\infty_\R(\D)$ is 2.
 \end{thm}
 \begin{proof}

 First we note that $\text{tsr}\,H^\infty_\R(\D)\not=1$, since, in view of Rouch\'e's theorem, 
  the function $z$ cannot be uniformly approximated by invertibles. 

 Now we will split our poof into several steps.  One step was contained in
Lemma \ref{sua}.  The idea is to approximate a function $\varphi\in H^\infty_\R(\D)$ by a product
of the form $\varphi_0=bFF^*$, 
where $b$ is a \CNBP\  with real zeros only
 and where $F\in H^\infty(\D)$ has no real zeros (see Lemma \ref{sua}).
  \medskip

 {\it Step 1.}~~ Let $g\in H^\infty_\R(\D)$ and let $p$ be  a \CNBP\
with real zeros.
 Then it can  be shown as in Lemma \ref{sua}, the first Step, that for every $\epsilon>0$ there exists
$G\in H^\infty_\R(\D)$
 such that $(G,p)\in U_2(H^\infty_\R(\D))$ and $||G-g||_\infty<\epsilon$.  For later purposes, it  is important that the function $p$ did not change!
\medskip

{\it Step 2.}~~  Let $(ff^*,gg^*)\in H^\infty_\R(\D)^2$ be a pair, both
entries $f$ and $g$ in
$H^\infty(\D)$ having no real zeros.
Since the topological stable rank of $H^\infty(\D)$ is two, see \cite{su},
for every $\epsilon\in\,]0,1[$
there exist $f_1$ and $g_1$ in $H^\infty(\D)$ such that
$$
(f_1,g_1)\in U_2(H^\infty(\D)), ~~ ||f-f_1||_\infty<\frac{\epsilon}{2}, ~~ ||g-g_1||_\infty<\frac{\epsilon}{2}.
$$
We may assume that $|f_1(z)|+|g_1(z)|\geq\delta>0$ for all $z\in\D$. Note that $\delta$
depends on $\epsilon$.

Consider now the pair $(f_1^*, g_1)$. By the same reason as above,
there exist $f_2$ and $g_2$ in $H^\infty(\D)$ such that $(f_2,g_2)\in U_2(H^\infty(\D))$,
$$
||f_1^*-f_2||_\infty<\min\left(\frac{\delta}{4},\frac{\epsilon}{2}\right), \textnormal{ and }
||g_1-g_2||_\infty<\min\left(\frac{\delta}{4},\frac{\epsilon}{2}\right).
$$

We claim that  $(f_2^*, g_2)\in U_2(H^\infty(\D))$, too.  In fact,
\begin{eqnarray*}
\abs{f_2^*(z)}+\abs{g_2(z)} & \geq & \abs{f_1(z)}+\abs{g_1(z)} - \bigl( ||f_1-f_2^*||_\infty + ||g_1-g_2||_\infty\bigr)\\
 & \geq & \delta-\frac{\delta}{2} =\frac{\delta}{2}.
\end{eqnarray*}

 Hence, $(f_2f_2^*, g_2)\in U_2(H^\infty(\D))$ and
 $$||f_2-f^*||_\infty \leq ||f_2-f_1^*||_\infty+||f_1^*-f^*||_\infty< \epsilon$$
 as well as
 $$||g_2-g||_\infty\leq ||g_2-g_1||_\infty +||g_1-g||_\infty<\epsilon.$$

 Noticing that $||f_2||_\infty\leq 1+||f_1||_\infty\leq 2+||f||_\infty$, we obtain that
 $$
 ||ff^*-f_2f_2^*||_\infty\leq ||(f^*-f_2)f||_\infty +||f_2(f-f_2^*)||_\infty$$
 $$
 \leq 2\epsilon(1+||f||_\infty).$$

 Since $f_2f_2^*\in H^\infty_\R(\D)$, we  get that
$(f_2f_2^*,g_2g_2^*)\in U_2(H^\infty_\R(\D))$.
 As above, we see that  $||g_2g_2^*-gg^*||_\infty\leq 2\epsilon(||g||_\infty+1)$.\medskip

 {\it Step 3.} ~~  Let $(ff^*,h)$ be a pair in $H^\infty_\R(\D)^2$, where
$f\in H^\infty(\D)$ has no real zeros. 
As mentioned at the beginning of this proof, we may
 approximate $h$ by  a product of the form $h=Cgg^*$, where
$C$ is a \CNBP\ with real zeros only and where $g\in H^\infty(\D)$
has no real zeros.

Continuing, it suffices then to consider  the pair $(ff^*,Cgg^*)$.
 Using Step 2, we find $F$ and $G \in H^\infty(\D)$ with $(FF^*,GG^*)\in U_2(H^\infty_\R(\D))$ and $||f-F||_\infty<\epsilon$ and $||g-G||_\infty<\epsilon$.

  We may assume that $|F(z)F^*(z)|+|G(z)G^*(z)|\geq\delta>0$
  for all $z\in \D$.
 Now we apply Step 1 to  the pair
 $(FF^*,C)$ and get  $R\in H^\infty_\R(\D)$ with $||R-FF^*||_\infty<
\min\left(\frac{\delta}{2}, \epsilon\right)$
 so that  $(R,C)\in U_2(H^\infty_\R(\D))$. We claim that $(R, CGG^*)\in
U_2(H^\infty_\R(\D))$,
 too. In fact, assuming that for some $x\in M(H^\infty(\D))$ we have
$R(x)=G(x)G^*(x)=0$,
 then $|F(x)F^*(x)|<\frac{\delta}{2}$ and so $|F(x)F^*(x)|+|G(x)G^*(x)|<\frac{\delta}{2}$, a
contradiction.

 Moreover,
 $$ ||R-ff^*||_\infty\leq ||R-FF^*||_\infty+||FF^*-ff^*||_\infty$$
 $$\leq ||R-FF^*||_\infty+||F^*(F-f)||_\infty+ ||f(F^*-f^*)||_\infty <  2\epsilon(1+||f||_\infty)$$
 and $$||CGG^*-Cgg^*||_\infty< 2\epsilon(1+||g||_\infty).$$

 Thus we have approximated $(ff^*,h)$ by an invertible pair $(R, CGG^*)$
in  $H^\infty_\R(\D)^2$.
 \medskip

 {\it Step 4.} ~~ Let us consider now an arbitrary pair $(\varphi, h)$
  in $H^\infty_\R(\D)^2$. As above, we approximate  $\varphi$ by the product
  $\varphi_0=Bff^*$, where $B$ is a \CNBP\ with real zeros only and where $f\in H^\infty(\D)$ has no real zeros.

 Continuing, it suffices to consider  the pair $(Bff^*,h)$.
 By Step 1 applied to $(B,h)$  there exists $H\in H^\infty_\R(\D)$ such
that
 $(B,H)\in U_2(H^\infty_\R(\D))$ and $||H-h||_\infty<\epsilon$. Assume that
$|B(z)|+|H(z)|\geq \delta$ for all $z\in\D$.
 Step 3 applied to $(ff^*,H)$ yields a pair $(R, K)\in
U_2(H^\infty_\R(\D)) $ such that
 $||R- ff^*||_\infty <\epsilon$ and $||H-K||_\infty<\min\left(\epsilon,\frac{\delta}{2}\right)$. Thus $(B,K)\in
U_2(H^\infty_\R(\D))$ as well,  since otherwise $B(x)=K(x)=0$ for some $x\in M(H^\infty(\D))$ would imply that
 $|B(x)|+|H(x)|<\frac{\delta}{2}$, a contradiction. Hence $(Bff^*, h)$, and so $(\varphi,h)$ has been approximated
 by the invertible pair $(BR,K)$ in $H^\infty_\R(\D)^2$.

 \end{proof}

 \begin{rem}
 Since the Bass stable rank is always less than or equal  the topological stable
rank, see \cite{CL},
 we get a second
 proof of Theorem \ref{str} above.  We have chosen to proceed as in the
present case, though,
 in order to make our proof of the stable rank determination of
$H^\infty_\R(\D)$ independent
 of the topological stable rank of $H^\infty_\R(\D)$.
 \end{rem}

\section{Bass and Topological Stable Ranks of $A_\R(\D)$}

We present the new proof of the results by Rupp and Sasane.

\begin{thm}[Rupp and Sasane, \cite{RS1}]

The Bass stable rank of $A_\R(\D)$ is 2.
\end{thm}

\begin{proof}

As above for the case of $H^\infty_\R(\D)$, we see that the stable rank of $A_\R(\D)$ is at least two.

Now, let $(f,g,h)$ be an invertible triple in $A_\R(\D)^3$. Since the stable rank of $H^\infty_\R(\D)$ is two, there exists $p,q,x,y\in H^\infty_\R(\D)$ such that 
$$1=x(f+ph)+y(g+qh).$$
 Now consider the dilations
$F_r$ given by $F_r(z):=F(rz)$, $0<r<1$, of all the seven functions appearing here. Since $f,g,h$ are in $A(\D)$, we may choose $r$ so close to 1 that
$$||f-f_r||_\infty<\frac{1}{5||x||_\infty+1},~~ ||g-g_r||_\infty<\frac{1}{5||y||_\infty+1}$$
and
$$
||h-h_r||_\infty<\min\left(\frac{1}{5||xp||_\infty+1},\frac{1}{5||yq||_\infty+1}\right).$$

Then $v:=x_r(f+p_rh)+y_r(g+q_rh)$ is invertible in $A_\R(\D)$, since $|v|>1/5$ on $\overline\D$.
\end{proof}

\begin{thm}[Rupp and Sasane, \cite{RS1}]

The topological stable rank of $A_\R(\D)$ is 2.
\end{thm}

\begin{proof}
We first observe that the topological stable rank of $A_\R(\D)$ is at least two,  since by Rouch\'e's Theorem, the function $z$ cannot be uniformly approximated by invertibles.  Next, let $(f,g)$ be a pair in $A_\R(\D)^2$.  Since the topological stable rank of $H^\infty_\R(\D)$ is two, for any $\epsilon>0$ there exists an invertible pair $(p,q)$ in $H^\infty_\R(\D)^2$ such that
$$
\norm{f-p}_\infty+\norm{g-q}_\infty<\frac{\epsilon}{2}.
$$
Again, for a function $F$, let $F_r(z):=F(rz)$, $0<r<1$, denote the dilation of the function.
Then, for every $0<r<1$
$$
\norm{f_r-p_r}_\infty+\norm{g_r-q_r}_\infty<\frac{\epsilon}{2}.
$$

 Note that  the functions $p_r$ and $q_r$ are in $A_\R(\D)$ and  that $(p_r,q_r)$ is an invertible pair in $A_\R(\D)^2$. 
 Since every disc algebra function is a uniform limit of its dilations,  we may choose $0<r_0<1$ so that 
$$
\norm{f_{r_0}-f}_\infty+\norm{g_{r_0}-g}_\infty<\frac{\epsilon}{2}.
$$
Then the pair $(p_{r_0},q_{r_0})$ is invertible in $A_\R(\D)^2$ and 
\begin{eqnarray*}
\norm{f-p_{r_0}}_\infty+\norm{g-q_{r_0}}_\infty & \leq & \norm{f-f_{r_0}}_\infty+\norm{g-g_{r_0}}_\infty\\ 
& + & \norm{f_{r_0}-p_{r_0}}_\infty+\norm{g_{r_0}-q_{r_0}}_\infty\\
 & < & \frac{\epsilon}{2}+\frac{\epsilon}{2}=\epsilon.
\end{eqnarray*}
Thus, the topological stable rank of $A_\R(\D)$ is two.
\end{proof}

Finally, we present an elementary  $\overline\partial$-free proof of the following result of the second author \cite{Wick2}. To this end, recall that  $U_2(A_\R(\D))$ denotes the set of all invertible pairs in $A_\R(\D)^2$.

\begin{thm} [Wick, \cite{Wick2}] 

An invertible pair $(f,g)$ in  $A_\R(\D)^2$ is reducible if and only if  $f$ has constant sign on the real zeros of $g$, or, equivalently, if and only if between two real zeros of $g$ in $\overline\D$ there is
an even number of zeros of $f$ (or none).
\end{thm}

\begin{proof}
Using a variant of the argument given in the introduction, we see that the condition on the zeros of $f$ and $g$ is necessary.  So, we turn to the sufficiency.

Let $f=BF$ be the Riesz factorization of $f$. Obviously, the zeros of $B$ are symmetric with respect to the real axis. Hence $F\in A_\R(\D)$, too.  Since $F$ has no zeros in $\D$, the dilates
$F_r$ given by $F_r(z)=F(rz)$ are invertible in $A_\R(\D)$; so the pairs  $(F_r, g)$ are reducible
in $A_\R(\D)$. Since the $F_r$ converge uniformly to $F$, we obtain from Lemma \ref{clopen}
that $(F,g)$ is reducible, too.  In the case where $B$ is an infinite \BP, let $B_n$ be a symmetric tail of the   $B$
(i.e. non-real zeros appear in pairs: $a$ and $\overline a$, so that $B_nF\in A_\R(\D)$.) Since $BF$ is in the disc algebra,
it is easy to see that $B_nF$ converges uniformly to $F$. Hence, for some $n_0$, we
conclude that for $n\geq n_0$, the pairs $(B_nF, g)$ are reducible in $A_\R(\D)$.

In the following, we let $B_{n_0}=1$ if $B$ is  a finite \BP.
 So let $(a_j, \overline a_j)$,   $(j=1,\dots,q)$,  be the pairs of non-real zeros of $\frac{B}{B_{n_0}}$, called the remaining zeros. 
Moreover, let $r_j, (j=1,\dots,m),$ be the remaining real zeros of $\frac{B}{B_{n_0}}$, ordered increasingly. 
Note that the hypothesis that $BF$ has constant sign on the real zeros of $g$ implies that
between two  consecutive real zeros  of $g$ there are an even number of zeros of $B$.

Since the reducibility of $(f_1,g)$ and $ (f_2,g)$ implies that $(f_1f_2,g)$ is reducible, too, 
in order to show that $(BF,g)$ is reducible,  it suffices to show  that for the remaining
  pairs of non-real zeros $(a_j, \overline a_j)$ of $B$ the pairs $\bigl((z-a_j)(z-\overline a_j), g\bigr)$
are reducible  and  that for the remaining real zeros $r_j$ the pair 
$\bigl((z-r_1)\cdots (z-r_{m}), g\bigr)$ is reducible.
Note that in the first  case, the abscissa $(x-a_j)(x-\overline a_j)= |x-a_j|^2 >0$ for $-1\leq x\leq 1$, hence
 are positive  on the real zeros of $g$. 
 
 Concerning the factor $(z-r_1)\cdots (z-r_{m})$, we also know that  between two consecutive zeros of $g$
 there is an even number (or none) of the $r_j$'s (since this is the case for $B$ and $B_n$). 
 
   We will use Lemma \ref{clopen} to show the reducibility of  $\bigl((z-a_j)(z-\overline a_j), g\bigr)$.
   
   Fix $j\in\{1,\dots, q\}$.
   Since $Z(g)$ is totally disconnected in $\D$, there is a Jordan  arc $\Gamma$ in $\overline \D$
   that avoids $Z(g)$ and satisfies $\Gamma(0)=a_j$ and $\Gamma (1)=e^{i\theta}$ for some
    $\theta\in [0,2\pi[$.  Consider the curve $\phi(t)=(z-\Gamma(t))(z-\overline{\Gamma(t)})$ inside $A_\R(\D)$.
    Note that for all $t$ we have that $(\phi(t),g)\in U_2(A_\R(\D))$. Moreover,
     $\phi(1)$ is an outer function in $A_\R(\D)$. By the first step of this proof, 
    the pair  $(\phi(1),g)$ 
     is reducible.  Hence by Lemma  \ref{clopen}, $$\bigl((z-a_j)(z-\overline a_j), g\bigr)=(\phi(0), g)$$ is reducible, too.

    The proof that $\bigl((z-r_1)\cdots(z-r_{m}), g\bigr)$ is reducible is  done similarly,
     but we need to split it into several cases:
        
    {\it Case 1}:  Suppose that between the consecutive zeros $w_1$ and $w_2$ of $g$ on $[-1,1]$
    there are $2p$ of the remaining  zeros of  $B$, denoted here by $r_1,\dots,r_{2p}$.
    
    Consider  a path $\Gamma$  in $\D$ that avoids $Z(g)$ and satisfies
    $\Gamma(0)=r_1$ and $\Gamma(1)=e^{i\theta}$ for some
    $\theta\in [0,2\pi[$. Then let  $\phi$ be  the curve 
    $$
    \phi(t)= (z-\Gamma(t))^p(z-\overline{\Gamma(t)})^p, \;0\leq t\leq 1.
    $$
    It is clear that $(\phi(t),g)\in U_2(A_\R(\D))$ for every $t$. Since $\phi(1)$ is an outer function,
    we  obtain, as above, that  $(\phi(1),g)$, and hence, via Lemma \ref{clopen}, that
    $(\phi(0), g)=\bigl((z-r_1)^{2p},g)$ is reducible.  (Note, that 
    at this point we have to use even powers).
    
    Now, consider   for $0\leq t\leq 1$  the curve

 $$
\phi(t)=\prod_{j=1}^{2p}(z-[r_1+t(r_j-r_1)])
 $$
 
 Again, $(\phi(t), g)$ is an invertible  pair in $A_\R(\D)^2$ (since there are no zeros of $g$ between
$r_1$ and $r_{2p}$.)  Since for $t=0$ we know that  $(\phi(0),g)=\bigl((z-r_1)^{2p},g\bigr)$
is reducible, we use Lemma  \ref{clopen} to deduce that 
$$ 
(\phi(1),g) = \bigl( (z-r_1)(z-r_2)\cdots(z-r_{2p}), g\bigr)
$$
is reducible. 
\medskip

{\it Case 2:} Suppose that $-1\leq w<1$ is  the biggest  zero of $g$ and that 
$w<r_ j< r_{j+1}<\dots<r_{m}$.

Now let us consider for $0\leq t\leq 1$ the curve

$$
\phi(t)=\prod_{k=j}^m(z-[r_k+t(1-r_k)])
$$

Then $\phi(1)$ is the  outer function $(z-1)^{m-j+1}$; hence $(\phi(1), g)$ is reducible.
Since $(\phi(t),g)\in U_2(A_\R(\D))$ for every $t$,
by Lemma \ref{clopen}, $(\phi(0),g)=\bigl( (z-r_j)(z-r_{j+1})\dots(z-r_m), g\bigr)$ is reducible, too.
\medskip

{\it Case 3:} Suppose that $-1<w\leq 1$ is the smallest zero of $g$ and that
$r_1<r_2<\dots<r_k<w$. Then we proceed as in the second case to prove
the reducibility of $\bigl( (z-r_1)(z-r_2)\dots (z-r_k), g\bigr)$.

\medskip

Now by multiplying all the factors of type $(z-r_l)\cdots(z-r_{l'})$ considered here 
in the three cases, we obtain
the reducibility of the pair $\bigl( (z-r_1)\cdots(z-r_m), g\bigr)$. 

\end{proof}

\section{Concluding Remarks}

While Theorem \ref{str} settles the conjecture of S. Treil regarding the Bass stable rank of $H^\infty_\R(\D)$, it does so in an unusual  manner.  Typically, when solving Bezout equations one would like to do so in such a way as to provide estimates of the solutions.  Thus, we conjecture that the following should be true,

\begin{conj}\label{conjec}
Let $f_1, f_2, f_3\in H^\infty_{\R}(\D)$ be such that $$1\geq\inf_{z\in\D}(\abs{f_1(z)}+\abs{f_2(z)}+\abs{f_3(z)})=\delta>0.$$
  Then there exists $h_1, h_2, g_1, g_2\in H^\infty_{\R}(\D)$ such that
$$
1=g_1(f_1+h_1f_3)+g_2(f_2+h_2f_3)
$$
and $\norm{g_j}_\infty,\norm{h_j}_\infty\leq C(\delta)$ for $j=1,2$.
\end{conj}

We note that there is a result of Gamelin \cite{gam} (see also \cite{Garnett}, p.~368) that 
the mere assumption that the unit disc is dense in the spectrum of  $H^\infty(\D)$  automatically
implies that every Corona problem in $\D$ admits a solution with norm control.  Can something similar
be said here? Note that the following observation is easy to show:

\begin{obs}
Let $\Omega=\bigcup_{n=0}^\infty (\D+ 3n)$ and let $H_\R(\Omega)$ be the set
of bounded analytic functions on $\Omega$ that are real on the reals within $\Omega$.
Then  Conjecture \ref{conjec} has a positive answer  if and only if every
such problem  has a bounded solution on $\Omega$ (without norm estimates).
\end{obs}

\end{document}